\documentclass{article}
\usepackage{amssymb,latexsym}
\usepackage{graphicx}
\usepackage{amsmath}
\usepackage{amsthm}
\usepackage{amscd}

\usepackage{bm}
\begin{document}
\newcommand{\D}{\displaystyle}
\newcommand{\DF}[2]{\frac{\D#1}{\D#2}}
\def\Xint#1{\mathchoice
{\XXint\displaystyle\textstyle{#1}}%
{\XXint\textstyle\scriptstyle{#1}}%
{\XXint\scriptstyle\scriptscriptstyle{#1}}%
{\XXint\scriptscriptstyle\scriptscriptstyle{#1}}%
\!\int}
\def\XXint#1#2#3{{\setbox0=\hbox{$#1{#2#3}{\int}$}
\vcenter{\hbox{$#2#3$}}\kern-.5\wd0}}
\def\ddashint{\Xint=}
\def\dashint{\Xint-}
\pagestyle{plain}
\newpage
\newtheorem{theorem}{Theorem}[section]
\newtheorem{definition}{Definition}[section]
\newtheorem{remark}{Remark}[section]
\newtheorem{corollary}{Corollary}[section]
\newtheorem{proposition}{Proposition}[section]
\newtheorem{lemma}{Lemma}[section]
\normalsize
\title{An Energy Reducing Flow for Multiple-Valued Functions\thanks{The author wants to thank his
advisor, Professor Robert Hardt for introducing this subject to him and numberous fruitful, enjoyful discussions.}}
\author{\large Wei Zhu\\
                          }
\maketitle
\begin{abstract}
By the method of discrete Morse flows, we construct an energy
reducing multiple-valued function flow. The flow we get is
H$\ddot{o}$lder continuous with respect to the $\mathcal{L}^2$
norm. We also give another way of constructing flows in some
special cases, where the flow we get behaves like ordinary heat
flow.
\end{abstract}
\newpage
\tableofcontents
\newpage
\section{Introduction}
This work was originally motivated by a paper $\cite{cx}$, in
which he constructed a mass reducing flow for integral currents.
The ideal of his construction comes from Horihata and Kikuchi's
paper $\cite{hk}$ for nonlinear parabolic equations.  More
specifically, let $T_0$ be an integral current, and let $h>0$ be
given. He defines a step approximation sequence $\{T_h^k\}$ of
integral currents with the same boundary as that of $T_0$ by
choosing $T_h^k$ such that $T_h^k$ minimizes the functional
$$G(T)=G(T_h^{k-1},T,h)=M(T)^2+\DF{F_k(T_h^{k-1}-T)^2}{h},$$
where $T$ is an integral current with $\partial T=\partial T_0$.
Then he constructed a $k+1$ dimensional current $S_h$ by
``connecting'' those $T_h^k$. Finally he takes a weak limit $S$ of
$S_h$ as $h\rightarrow 0$, slices $S$ with respect to $t$ to get
an integral current at time $t$. The flow is H$\ddot{o}$lder
continuous under the flat norm and reduces the mass of the initial
integral current while keeping the boundary fixed. Later on, this
same time discretization process was used in Haga, Hoshino and
Kikuchi's paper $\cite{hhk}$ to construct a harmonic map flow.
This gives an alternative proof of the classical result due to J.
Eells, Jr. and
J.H.Sampson $\cite{es}$ and to R.S.Hamilton $\cite{hrs}$.\\
Our work is trying to construct a what is so called
``multiple-valued harmonic flow" by the similar time
discretization method. One obvious obstacle here is that we do not
have differential equations for multiple-valued functions.
Therefore a lot of PDE methods can not be applied. Another thing
that stands in the way is due to the lack of fundamental algebraic
operations, for example, addition for multiple-valued functions.
Hence we can not use linear interpolation to connect those
multiple-valued function at different stages. All of them will be handled with great care.\\
In the first part of this paper, a basic description of
multiple-valued functions and related results are given. Most of
them can be found in $\cite{af}$.\\
In the second part, we state and prove some theorems that we need
for
the construction later.\\
In the third part, an energy reducing flow of multiple-valued
functions is constructed. Those properties of this flow are
proved.\\
Finally, we look at the special case when
$\mathbb{Q}=\mathbb{Q}_2(\mathbb{R})$ and the initial
multiple-valued function $f_0$ satisfies $\eta(f_0)\equiv 0$, i.e
symmetric. We will show that the flow will be symmetric all the
time and each component separates immediately. A concrete example
when the domain is two dimensional is also given.
\section{Preliminaries}
The theory of multiple-valued functions was developed in
$\cite{af}$. It is the most natural framework for the
regularity theory in geometric measure theory and promises a lot of future development and applications in other fields. Here we introduce the
basic notations and facts of multiple-valued functions. The readers are referred to $\cite{af}$ for more details. We also use standard terminology in geometric measure theory,
all of which can be found on page 669-671 of the treatise {\it Geometric Measure Theory} by H. Federer $\cite{fh}$.\\
\newline
The space $\mathbb{Q}=\mathbb{Q}_Q(\mathbb{R}^n)$ consists of all the unordered $Q$ points in $\mathbb{R}^n$, denoted by
$\sum_{i=1}^Q [[p_i]],p_i\in \mathbb{R}^n$. We let $\mbox{spt}(\sum_{i=1}^Q [[p_i]])=\cup_{i=1}^Q \{p_i\}$.\\
\newline
We define a metric on $\mathbb{Q}$
$$\mathcal{G}:\mathbb{Q}\rightarrow \mathbb{R}$$
by setting for $p_1,\cdot\cdot\cdot,p_Q,q_1,\cdot\cdot\cdot,q_Q\in\mathbb{R}^n$,
$$\mathcal{G}(\sum_i[[p_i]],\sum_i[[q_i]])=\inf_{\sigma}\{(\sum_i|p_i-q_{\sigma(i)}|^2)^{1/2}:\sigma\;\mbox{is a permutation
of}\;\{1,\cdot\cdot\cdot,Q\}\}.$$
We let $|\sum_{i=1}^Q [[p_i]]|^2=\sum_{i=1}^Q|p_i|^2=\mathcal{G}^2(\sum_{i=1}^Q[[p_i]],Q[[0]]).$\\
$$\zeta:\mathbb{O}^*(n,1)\times \mathbb{Q}\rightarrow \mathbb{R}^Q\cap\{s:s_1\le s_2\le s_3\le\cdot\cdot\cdot\le s_Q\}$$
is defined for $\pi\in\mathbb{O}^*(n,1),p\in\mathbb{Q}$ by requiring $-\infty<s_1\le s_2\le s_3\le\cdot\cdot\cdot\le s_Q<\infty$
and $\pi_\sharp p=\sum_{i=1}^Q[[s_i]]$. One can easily check that Lip$(\zeta(\pi,\cdot))$=1 for each $\pi\in\mathbb{O}^*(n,1)$.\\
\newline
$P$ is the positive integer characterized in the following manner: if
$$b=\inf\{\arctan(1/n!),(\sin[2^{-1}\arctan(1/n!)])^{n-1}/2QQ!(n-1)\},$$
and $c$ denotes the unique positive integer for which
$$1/b\le c<1+1/b,$$
then
$$P=2^{-1}([4c(n-1)+1]^n-[4c(n-1)-1]^n).$$
\begin{theorem}[$\cite{af},\S 1.2$]
There exist $\Pi_1,\Pi_2,\cdot\cdot\cdot,\Pi_p\in\mathbb{O}^*(n,1)$ such that\\
(1) $\Pi_i(x)=x_i$ for each $i=1,\cdot\cdot\cdot,n$ and each $x\in\mathbb{R}^n$.\\
(2) Lip$(\xi_0)=1$.\\
(3) $\xi:\mathbb{Q}\rightarrow\mathbb{Q}^*$ is a bilipschitzian homeomorphism with Lip$(\xi)\le P^{1/2}$ and
Lip$(\xi^{-1}|\mathbb{Q}^*)\le 1/(n^{1/2}\sin(b/2))$ corresponding to $b$ as above, where\\
$$\xi=\zeta(\Pi_1,\cdot)\Join\cdot\cdot\cdot\Join\zeta(\Pi_p,\cdot):\mathbb{Q}\rightarrow \mathbb{Q}^{PQ},$$
$$\mathbb{Q}^*=\xi(\mathbb{Q})$$
and
$$\xi_0=\zeta(\Pi_1,\cdot)\Join\cdot\cdot\cdot\Join\zeta(\Pi_n,\cdot):\mathbb{Q}\rightarrow \mathbb{R}^{nQ}.$$
\end{theorem}
(Here we use the notation that whenever $f:A\rightarrow B$ and
$g:A\rightarrow C$, we define
$$f\Join g:A\rightarrow B\times C, (f\Join g)(a)=(f(a),g(a)), a\in A.)$$
\begin{remark}Brian White showed that there is a modified bilipschitzian
correspondance
$\xi:\mathbb{Q}\rightarrow\mathbb{Q}^*\subset\mathbb{R}^{PQ}$ such
that for every $p\in\mathbb{Q}$, $p$ has a small neighbourhood in
$\mathbb{Q}$ such that $\xi$ is an equidistance map over the
neighbourhood. The modification is to choose the orthogonal
projections $\Pi_1,\cdot \cdot\cdot,\Pi_P$ in $\cite{af}$ as
complete sets of coordinate projections corresponding to distinct
orthonormal coordinate systems for $\mathbb{R}^n$ and to compose
the resulting map $\xi$ there with proper scaling to get such a
$\xi$. It has some other useful properties that we will mention
later. Moreover, we will use the modified $\xi$ throughout the
rest of this paper.
\end{remark}
\begin{theorem}[$\cite{af},\S 1.2$]
Suppose $-\infty<r(1)\le r(2)\le \cdot\cdot\cdot\le r(Q)<\infty$
and $-\infty<s(1)\le s(2)\le\cdot\cdot\cdot\le s(Q)<\infty$. Then
$$\sum_{i=1}^Q [r(i)-s(i)]^2=\inf_{\sigma}\{\sum_{i=1}^Q [r(i)-s(\sigma(i))]^2:\;\sigma\;\mbox{is a permutation
of}\;\{1,\cdot\cdot\cdot,Q\}\}.$$
\end{theorem}
\begin{remark}
This theorem says that the distance between two elements in
$\mathbb{Q}_Q(\mathbb{R})$ is obtained by matching those Q-tuples
pairwise according to the ascending order.
\end{remark}
\begin{theorem}[$\cite{af},\S 1.3$]
There exists an explicitly constructable, piecewise linear function
$$\rho:\mathbb{R}^{PQ}\rightarrow \mathbb{R}^{PQ}$$
such that Lip$(\rho)<\infty$, $\rho(\mathbb{R}^{PQ})\subset\mathbb{Q}^*$, and $\rho(x)=x$ for each $x\in\mathbb{Q}^*$.
\end{theorem}
\begin{definition} (a) $f$ is a Q-valued function on some subset
$U$ of $\mathbb{R}^m$ if it is a map
$$f:U\subset \mathbb{R}^m\rightarrow \mathbb{Q}$$
(b) For a given smooth, compact embedded manifold $N$ in $\mathbb{R}^n$,
$$\underline{Q}(N)=\{\sum_{i=1}^Q[[p_i]],p_i\in N,i=1,\cdot\cdot\cdot,Q\}$$
(c) $f$ is a Q-valued map from some subset $U$ of $\mathbb{R}^m$
into $N$ if it is a map
$$f:U\subset \mathbb{R}^m\rightarrow \underline{Q}(N)\subset \mathbb{Q}$$
(d) Similarly we can define $\underline{Q}(V)$ for any vector
space $V$.
\end{definition}
\begin{definition}(a) $f$ is called a Q-valued Lipschitz function(map) if there
is a constant $C>0$ such that
$$\mathcal{G}(f(x),f(y))\le C|x-y|,\;\;x,y\in U.$$
(b) $f$ is called affine if there are $A_1,\cdot\cdot\cdot,A_Q$ where each $A_i$ is an affine map from
$\mathbb{R}^m$ to $\mathbb{R}^n$, such that
$$f(x)=\sum_{i=1}^Q [[A_i(x)]].$$
(c) $f$ is called affinely approximatable at $x_0$ if there are affine maps $A_1,\cdot\cdot\cdot,A_Q$ from $\mathbb{R}^m$
to $\mathbb{R}^n$ such that
$$\lim_{|x-x_0|\to 0}\DF{\mathcal{G}(f(x),\sum_{i=1}^Q [[A_i(x)]])}{|x-x_0|}=0.$$
(d) $f$ is strongly affinely approximatable at $x_0$ if (c) holds
for $f$ at $x_0$ and $A_i=A_j$ if $A_i(x_0)=A_j(x_0)$.
\end{definition}
\begin{remark}(1) From $\cite{af},\S 1.4$, if $f$ is a Q-valued Lipschitz function, then it is strongly affinely approximatable almost everywhere over its domain.\\
(2) If $f$ is affinely approximatable at $x_0$ with
$\sum_{i=1}^Q[[A_i]]$ as its affine approximation, then obviously
$f(x_0)= \sum_{i=1}^Q [[A_i(x_0)]]$ and
$A_i(x)=A_i(x_0)+L_i(x-x_0)$ with
$L_i\in\;\mbox{Hom}(\mathbb{R}^m,\mathbb{R}^n)$.
\end{remark}
\begin{definition} If $f$ is affinely approximatable at $x_0$, then\\
(a) $\sum_{i=1}^Q[[L_i]]\in
\underline{Q}(\mbox{Hom}(\mathbb{R}^m,\mathbb{R}^n))$, denoted by
$Df(x_0)$ is defined as the differential of $f$ at $x_0$. We let
$|Df(x_0)|^2=\sum_{i=1}^Q|L_i|^2$, where $|L|$ is the Euclidean
norm of the
matrix associated with any $L\in\;\mbox{Hom}(\mathbb{R}^m,\mathbb{R}^n).$\\
(b) $\sum_{i=1}^Q[[L_i(v)]]$ is defined as the derivative of $f$
at $x_0$ in the direction $v$ and is denoted by $D_v f\in
\mathbb{Q}$. Let $|D_v f(x_0)|^2=\sum_{i=1}^Q |L_i(v)|^2$.
\end{definition}
\begin{remark}The map $\xi$ mentioned before has the following properties
$$|\xi\circ f|=|f|,|D_v(\xi\circ f)|=|D_v f|$$
\end{remark}
\begin{definition} (a) Suppose $A\subset \mathbb{R}^m$ is bounded and open and that
$\partial A$ is an $m-1$ dimensional submanifold of $\mathbb{R}^m$
of class 1. Whenever $V$ is a Euclidean vector space,
$\mathcal{Y}_2(\mathbb{R}^m,V)$ and $\mathcal{Y}_2(A,V)$ are the
real vector spaces of square summable functions whose distribution
first derivatives are also square summable. $\partial\mathcal{Y}_2
(\partial A,V)$ is the real vector space of all
$\mathcal{H}^{m-1}$ measurable function $f:\partial A\rightarrow
V$ such that
$$\int_{\partial A}|f|^2 d\mathcal{H}^{m-1}+\int_{z\in\partial A} |z|^{-m}\int_{x\in\partial A}|f(x+z)-f(z)|^2 d
\mathcal{H}^{m-1}x d\mathcal{H}^{m-1}z<\infty$$ (b) Whenever
$K\subset \mathbb{R}^m$ is $\mathcal{L}^m$ measurable with
$\mathcal{L}^m(K\sim A)=0$ and $f\in \mathcal{Y}_2(A,V)$ we set
$$\mbox{Dir}(f;K)=\int_K|Df|^2 d\mathcal{L}^m.$$
Additionally, we define for each $g\in\partial\mathcal{Y}_2(\partial A,V)$ and each $\mathcal{H}^{m-1}$ measurable
set $L\subset\mathbb{R}^m$ with $\mathcal{H}^{m-1}(L\sim \partial A)=0$,
$$\mbox{dir}(g;L)=\int_L |Dg|^2 d\mathcal{H}^{m-1}.$$
\end{definition}
\begin{definition}
Assuming $f:\mathbb{R}^m\rightarrow V$ is locally $\mathcal{L}^m$
summable, we say $f$ is strictly defined if and only if whenever
$x\in\mathbb{R}^m$, and there is some $y\in V$ for which
$$\lim_{r\downarrow 0}\dashint_{z\in B^m(x,r)}|f(z)-y|d\mathcal{L}^m z=0,$$
then $f(x)=y$. \end{definition}
\begin{definition}
Whenever $f\in\mathcal{Y}_2(A,V)$ and $g:\partial A\rightarrow V$
we say that $f$ has boundary values $g$ if and only if there is
$h\in\mathcal{Y}_2(\mathbb{R}^m,V)$ which is strictly defined such
that
$$\mathcal{L}^m(A\cap\{x:f(x)\not=h(x)\})=0=\mathcal{H}^{m-1}(\partial A\cap\{x:g(x)\not=h(x)\}).$$
\end{definition}
\begin{definition}
(a) We define
$$\mathcal{Y}_2(\mathbb{R}^m,\mathbb{Q})\;[\mbox{resp.}\;\mathcal{Y}_2(A,\mathbb{Q})]$$
to be the space of all functions $f:\mathbb{R}^m\rightarrow \mathbb{Q}$[resp. $f:A\rightarrow \mathbb{Q}$] such that
$\xi\circ f\in\mathcal{Y}_2(\mathbb{R}^m,\mathbb{R}^{PQ})$[resp. $\xi\circ f\in\mathcal{Y}_2(A,\mathbb{R}^{PQ})$]. We also
define
$$\partial\mathcal{Y}_2(\partial A,\mathbb{Q})$$
as the space of all functions $g:\partial A\rightarrow \mathbb{Q}$ such that $\xi\circ g\in\partial\mathcal{Y}_2(\partial A,
\mathbb{R}^{PQ})$.\\
(b) For each $f\in\mathcal{Y}_2(A,\mathbb{Q})$ and each
$\mathcal{L}^m$ measurable set $K\subset \mathbb{R}^m$ which is
$\mathcal{L}^m$ almost a subset of $A$, we define
$$Dir(f;K)=Dir(\xi_0\circ f;K).$$
For each $g\in\partial\mathcal{Y}_2(\partial A,\mathbb{Q})$ and
$\mathcal{H}^{m-1}$ measurable set $L\subset\mathbb{R}^m$ which is
$\mathcal{H}^{m-1}$ almost a subset of $\partial A$, we set
$$dir(g;L)=dir(\xi_0\circ g;L).$$ (c) Whenever
$f\in\mathcal{Y}_2(\mathbb{R}^m,\mathbb{Q})$[resp.
$f\in\mathcal{Y}_2(A,\mathbb{Q})$] we say that $f$ is strictly
defined if and only if $\xi\circ f$ is strictly defined. One notes
that in case $f:\mathbb{R}^m\rightarrow \mathbb{R}^{PQ}$ is
locally $\mathcal{L}^m$ summable with im$(f)\subset \mathbb{Q}^*$,
$x\in \mathbb{R}^m$, $y\in\mathbb{R}^{PQ}$, and
$$\lim_{r\downarrow 0}\dashint_{z\in B^m(x,r)}|f(z)-y|d\mathcal{L}^m z=0,$$
then $y\in\mathbb{Q}^*$ since $\mathbb{Q}^*$ is closed.\\ (d) For
$f\in\mathcal{Y}_2(A,\mathbb{Q})$,
$g\in\partial\mathcal{Y}_2(\partial A,\mathbb{Q})$ one says that
$f$ has boundary values $g$ if and only if $\xi\circ f$ has
boundary values $\xi\circ g$. \end{definition}
\begin{definition}One says that $f:A\rightarrow \mathbb{Q}$ is Dir
minimizing if and only if $f\in\mathcal{Y}_2(A,\mathbb{Q})$ and,
assuming $f$ has boundary values
$g\in\partial\mathcal{Y}_2(\partial A,\mathbb{Q})$, one has
$$\mbox{Dir}(f;A)=\inf\{\mbox{Dir}(h;A):h\in\mathcal{Y}_2(A,\mathbb{Q})\;\mbox{has boundary values}\;g\}.$$
\end{definition}
\begin{remark}All those definitions are easily extended to multiple-valued
maps. \end{remark}
\begin{theorem}[$\cite{af}, \S 2.2$]
Suppose $A\subset \mathbb{R}^m$ is bounded and open and $\partial
A$ is an $m-1$ dimensional submanifold of $\mathbb{R}^m$ of class
1. \\
(1) Let $f\in\mathcal{Y}_2(A,\mathbb{Q})$. Then\\
(a) for $\mathcal{L}^m$ almost all $x\in A$, $apD(\xi\circ f)(x)$
exists.\\
(b) if $x\in A$ and $apD(\xi\circ f)(x)$ exists, then both
$apD(\xi_0\circ f)(x)$ and $apAf(x)$ exist with $|apD(\xi_0\circ
f)(x)|=|apAf(x)|.$\\
(c) whenever $K\subset \mathbb{R}^m$ is $\mathcal{L}^m$ measurable
and is $\mathcal{L}^m$ almost a subset of $A$, $Dir(f;K)=\int_A
|apAf|^2 d\mathcal{L}^m$.\\ (2) Let $g\in\partial
\mathcal{Y}_2(\partial A,\mathbb{Q}).$ Then there exists
$f\in\mathcal{Y}_2(A,\mathbb{Q})
$ such that\\
(a) $f$ has boundary values $g$.\\
(b)
$\mbox{Dir}(f;A)=\inf\{\mbox{Dir}(h;A):h\in\mathcal{Y}_2(A,\mathbb{Q})\;\mbox{has
boundary values}\;g\}.$
\end{theorem}
\begin{theorem}[$\cite{af},\S A.1.2$]
Suppose $V$ is a Euclidean vector space and $A\subset \mathbb{R}^m$
is bounded and open such that $\partial A$ is a compact $m-1$
dimensional submanifold of $\mathbb{R}^m$ of class 1. Suppose $K$ is
a closed subset of $V$, $g\in\partial \mathcal{Y}_2(A,V)$,
$f_1,f_2,f_3,\cdot\cdot\cdot\in\mathcal{Y}_2(A,V)$ such that the following conditions hold:\\
(a) $g(x)\in K$ for each $x\in \partial A,$\\
(b) $f_i(x)\in K$ for each $x\in A$ and each $i=1,2,3,\cdot\cdot\cdot,$\\
(c) $f_i$ has boundary values $g$ for each $i=1,2,3,\cdot\cdot\cdot,$\\
(d) $\sup_i Dir(f_i;A)<\infty.$\\
Then there exists a subsequence $i_1,i_2,i_3,\cdot\cdot\cdot$ of
$1,2,3,\cdot\cdot\cdot$ and $f\in\mathcal{Y}_2(A,V)$ with the
following
properties:\\
(1) $0=\lim_{k\rightarrow \infty}\int_A|f-f_{i_k}|^2 d\mathcal{L}^m,$\\
(2) $\lim_{k\rightarrow \infty}\int_{x\in A}<ap Df_{i_k}(x)-ap Df(x),\phi(x)>d\mathcal{L}^m x=0\in V$ for each $\phi\in \mathcal{D}(A,\mathbb{R}^m)$,\\
(3) $f(x)\in K$ for each $x\in A$,\\
(4) $f$ has boundary values $g$,\\
(5) $\lim_{k\rightarrow
\infty}Dir(f_{i_k};A)=Dir(f;A)+\lim_{k\rightarrow
\infty}Dir(f_{i_k}-f;A)\in\mathbb{R}.$
\end{theorem}
\begin{remark}
It is easy to see that this theorem can be easily extended to the
case when $A$ is a cylinder of the form $[0,\infty)\times B_1^m$,
where $B_1^m=\{x:x\in\mathbb{R}^m, |x|\le 1\}$.
\end{remark}
\begin{definition}Define
$$(+):\mathcal{Y}_2(A,\mathbb{Q})\times
\mathcal{Y}_2(A,\mathbb{R}^n)\rightarrow
\mathcal{Y}_2(A,\mathbb{Q})$$ by setting
$$f(+)\phi(x)=(+)(f,\phi)(x)=\tau(-\phi(x))_\sharp f(x)$$
for each
$f\in\mathcal{Y}_2(A,\mathbb{Q})$,$\phi\in\mathcal{Y}_2(A,\mathbb{R}^n),x\in
A$. Here $\tau$ is the translation operator
$$\tau(y):\mathbb{R}^n\rightarrow \mathbb{R}^n,
\tau(y)(x)=x-y,\;\mbox{for}\;x\in\mathbb{R}^n.$$
\end{definition}
\begin{theorem}[$\cite{af},\S 2.3$]
$$Dir(f(+)\phi;A)=Dir(f;A)+2Q\int_A<D(\eta\circ
f),D\phi>d\mathcal{L}^m+QDir(\phi;A)$$ whenever
$f\in\mathcal{Y}_2(A,\mathbb{Q}),\phi\in\mathcal{Y}_2(A,\mathbb{R}^n)$.
\end{theorem}
\section{Extension of Luckhaus Lemma to Multiple-Valued Maps}
\begin{theorem}[Luckhaus Lemma, $\cite{sl}, \S 2.6$]
Suppose $N$ is an arbitrary compact subset of $\mathbb{R}^n$, $m\ge 2$, and $u,v\in W^{1,2}(\mathbb{S}^{m-1};N)$. Then there
is a constant $C$ which depends only on $m,n$ and $N$ such that
for each $\epsilon\in(0,1)$ there is a function $w\in W^{1,2}(\mathbb{S}^{m-1}\times[0,\epsilon];\mathbb{R}^n)$ such that $w$ agrees with $u$
in a neighbourhood of $\mathbb{S}^{m-1}\times\{0\}$, $w$ agrees with $v$ in a neighbourhood of $\mathbb{S}^{m-1}\times\{\epsilon\}$,
$$\int_{\mathbb{S}^{m-1}\times[0,\epsilon]}|\overline{D} w|^2\le C\epsilon\int_{\mathbb{S}^{m-1}}(|Du|^2+|
Dv|^2)+C\epsilon^{-1}\int_{\mathbb{S}^{m-1}}|u-v|^2,$$
and\\
$\mbox{dist}^2(w(x,s),N)\le$
$$C\epsilon^{1-m}(\int_{\mathbb{S}^{m-1}}|Du|^2+|Dv|^2)^{1/2}(\int_{\mathbb{S}^{m-1}}|u-v|^2)^{1/2}+C\epsilon^{-m}
\int_{\mathbb{S}^{m-1}}|u-v|^2$$
for $a.e.(x,s)\in \mathbb{S}^{m-1}\times[0,\epsilon].$ Here $D$ is the gradient on $\mathbb{S}^{m-1}$ and $\overline{D}$ is
the gradient on the product space $\mathbb{S}^{m-1}\times [0,\epsilon]$.
\end{theorem}
\begin{theorem}
Suppose $N$ is an arbitrary smooth compact manifold embedded in
$\mathbb{R}^n$, $m\ge 2$, and $u,v\in
\mathcal{Y}_2(\mathbb{S}^{m-1},\underline{Q}(N))$. Then there is a
constant $C$ which depends only on $m,n,Q$ and $N$ such that for
each $\epsilon\in(0,1)$ there is a map $w\in
\mathcal{Y}_2(\mathbb{S}^{m-1} \times[0,\epsilon];\mathbb{Q})$
such that $w$ agrees with $u$ in a neighbourhood of
$\mathbb{S}^{m-1}\times\{0\}$, $w$ agrees with $v$ in a
neighbourhood of $\mathbb{S}^{m-1}\times \{\epsilon\}$,
$$\int_{\mathbb{S}^{m-1}\times[0,\epsilon]}|\overline{D}w|^2\le C\epsilon\int_{\mathbb{S}^{m-1}}(|Du|^2+|
Dv|^2)+C\epsilon^{-1}\int_{\mathbb{S}^{m-1}}\mathcal{G}^2(u,v),$$
and\\
$\mbox{dist}^2(w(x,s),\underline{Q}(N))\le$
$$C\epsilon^{1-m}(\int_{\mathbb{S}^{m-1}}|Du|^2+|Dv|^2)^{1/2}(\int_{\mathbb{S}^{m-1}}\mathcal{G}^2(u,v))^{1/2}
+C\epsilon^{-m}
\int_{\mathbb{S}^{m-1}}\mathcal{G}^2(u,v)$$
for $a.e.(x,s)\in \mathbb{S}^{m-1}\times[0,\epsilon].$ Here $D$ is the gradient on $\mathbb{S}^{m-1}$ and $\overline{D}$ is
the gradient on the product space $\mathbb{S}^{m-1}\times [0,\epsilon]$.
\end{theorem}
\begin{proof}
Apply Luckhaus lemma to the set $\xi\circ \underline{Q}(N)$ and functions $\xi\circ u$, $\xi\circ v$, we get a function
$\tilde{w}\in W^{1,2}(\mathbb{S}^{m-1}\times[0,\epsilon];\mathbb{R}^{PQ})$ such that $\tilde{w}$ agrees with $\xi\circ u$ in a neighbourhood
of $\mathbb{S}^{m-1}\times\{0\}$, $\tilde{w}$ agrees with $\xi\circ v$ in a neighbourhood of $\mathbb{S}^{m-1}\times\{\epsilon\}$,\\
$$\int_{\mathbb{S}^{m-1}\times[0,\epsilon]}|\overline{D} \tilde{w}|^2 \le \tilde{C}\epsilon\int_{\mathbb{S}^{m-1}}(|D(\xi\circ u)|^2+|
D(\xi\circ v)|^2)+\tilde{C}\epsilon^{-1}\int_{\mathbb{S}^{m-1}}|\xi\circ u-\xi\circ v|^2,$$
and $\mbox{dist}^2(\tilde{w}(x,s),\xi\circ \underline{Q}(N))\le$
$$\tilde{C}\epsilon^{1-m}(\int_{\mathbb{S}^{m-1}}|D(\xi\circ u)|^2+|D(\xi\circ v)|^2)^{1/2}(\int_{\mathbb{S}^{m-1}}
|\xi\circ u-\xi\circ v|^2)^{1/2}$$
$$+\tilde{C}\epsilon^{-m}
\int_{\mathbb{S}^{m-1}}|\xi\circ u-\xi\circ v|^2$$
for $a.e.(x,s)\in \mathbb{S}^{m-1}\times[0,\epsilon].$ Here $D$ is the gradient on $\mathbb{S}^{m-1}$ and $\overline{D}$ is
the gradient on the product space $\mathbb{S}^{m-1}\times [0,\epsilon]$.\\
Now we define $w=\xi^{-1}\circ \rho\circ \tilde{w}\in
\mathcal{Y}_2(\mathbb{S}^{m-1}\times [0,\epsilon];\mathbb{Q})$. It
is easy to see that
$w$ agrees with $u$ in a neighbourhood of $\mathbb{S}^{m-1}\times\{0\}$ and $w$ agrees with $v$ in a neighbourhood of $\mathbb{S}^{m-1}\times\{\epsilon\}$.\\
The rest of the proof is obviously easy once we notice that
$\mbox{Lip}\;\xi^{-1}, \mbox{Lip} \;\xi, \mbox{Lip}\;\rho$ are all
finite, depending only on $m,n,Q$.
\end{proof}
\begin{corollary}
Suppose $N$ is a smooth compact manifold embedded in
$\mathbb{R}^n$, and $\Lambda>0$. There are $\delta_0
=\delta_0(m,n,N,Q,\Lambda)$ and $C=C(m,n,N,Q,\Lambda)$ such that the following hold:\\
(1) If we have $\epsilon\in(0,1)$ and if $u\in
\mathcal{Y}_2(B_\rho^m(y);\underline{Q}(N))$ with
$\rho^{2-m}\int_{B_\rho^m(y)}|Du|^2\le \Lambda$, and
$\epsilon^{-2m}\rho^{-m}\int_{B_\rho^m(y)}|\xi\circ
u-\lambda_{y,\rho}|^2\le \delta_0^2$, (where
$\lambda_{y,\rho}=\dashint_{ B_\rho^m(y)} \xi\circ u$) then there
is $\sigma\in (\DF{3\rho}{4},\rho)$ such that there is a map
$w=w_\epsilon\in \mathcal{Y}_2(B_\rho^m(y);\underline{Q}(N))$
which agrees with $u$ in a neighbourhood of $\partial
B_\sigma^m(y)$ and which satisfies
$$\sigma^{2-m}\int_{B_\sigma^m(y)}|Dw|^2\le \epsilon\rho^{2-m}\int_{B_\rho^m(y)}|Du|^2+\epsilon^{-1}C\rho^{-m}\int_{B_\rho^m(y)}
|\xi\circ u-\lambda_{y,\rho}|^2.$$ (2) If
$\epsilon\in(0,\delta_0]$, and if $u,v\in
\mathcal{Y}_2(B_{(1+\epsilon)\rho}^m(y)\backslash
B_\rho^m(y);\underline{Q}(N))$ satisfy the inequalities
$\rho^{2-m}\int_{B_{\rho(1+\epsilon)}^m(y)\backslash B_\rho^m(y)}
(|Du|^2+|Dv|^2)\le \Lambda$ and $\epsilon^{-2m}
\rho^{-m}\int_{B_{\rho(1+\epsilon)}^m(y)\backslash B_\rho^m(y)} \mathcal{G}(u,v)^2$\\
$<\delta_0^2$, then there is $w\in \mathcal{Y}_2(
B_{\rho(1+\epsilon)}^m(y)\backslash B_\rho^m(y);\underline{Q}(N))$
such that $w=u$ in a neighbourhood of $\partial B_\rho^m(y)$,
$w=v$
in a neighbourhood of $\partial B_{(1+\epsilon)\rho}^m(y)$ and\\
$\rho^{2-m}\int_{B_{\rho(1+\epsilon)}^m(y)\backslash B_\rho^m(y)}
|Dw|^2$
$$\le C\rho^{2-m}\int_{B_{\rho(1+\epsilon)}^m(y)\backslash B_\rho^m(y)} (|Du|^2+|Dv|^2)+C\epsilon^{-2}\rho^{-m}
\int_{B_{\rho(1+\epsilon)}^m(y)\backslash B_\rho^m(y)}
\mathcal{G}(u,v)^2.$$
\end{corollary}
\begin{proof}
The proof is basically the same as in $\cite{sl},\S 2.7$.
\end{proof}
\section{Compactness Theorem for Multi-Valued Minimizing Maps}
\begin{theorem}[Rellich Compactness Lemma, $\cite{sl}, \S 1.3$]
Suppose $\Omega$ is a bounded Lipschitz domain in $\mathbb{R}^m$ and $u_k$ is a sequence of $W^{1,2}(
\Omega)$ with \\$\mbox{sup}_k||u_k||_{W^{1,2}(\Omega)}<\infty$. Then there is a subsequence $u_k'$ and $u\in W^{1,2}(\Omega)$ such that
\begin{equation*}
\begin{split}
(a)&\;u_k'\rightharpoonup u\;\mbox{weakly in}\; W^{1,2}(\Omega),\\
(b)&\;u_k'\rightarrow u\;\mbox{strongly in}\; \mathcal{L}^2(\Omega),\\
(c)&\;\int_{\Omega} |Du|^2\le \liminf_{k'\to\infty}\int_{\Omega}|Du_k'|^2.
\end{split}
\end{equation*}
\end{theorem}
\begin{definition}Given a sequence $u_k$, and
$u\in\mathcal{Y}_2(\Omega,\underline{Q}(N))$, we say that,
$$u_k\rightharpoonup u\;\mbox{weakly in}\;\mathcal{Y}_2(\Omega)\;\mbox{if and only if}\;\xi\circ u_k\rightharpoonup
\xi\circ u\;\mbox{weakly in}\; W^{1,2}(\Omega)$$
$$u_k\rightarrow u\;\mbox{strongly in}\;L^2(\Omega)\;\mbox{if and only if}\;\xi\circ u_k\rightarrow \xi\circ u\;\mbox{strongly in}\;
\mathcal{L}^2(\Omega).$$
\end{definition}
\begin{theorem}
Suppose $\Omega$ is a bounded Lipschitz domain in $\mathbb{R}^m$ and $u_k$ is a sequence of $\mathcal{Y}_2(
\Omega,\mathbb{Q})$ with $\mbox{sup}_k||u_k||_{\mathcal{Y}_2(\Omega)}<\infty$. Then there is a subsequence $u_k'$ and $u\in \mathcal{Y}_2(\Omega,
\mathbb{Q})$ such that
\begin{equation*}
\begin{split}
(a)&\;u_k'\rightharpoonup u\;\mbox{weakly in}\; \mathcal{Y}_2(\Omega),\\
(b)&\;u_k'\rightarrow u\;\mbox{strongly in}\; \mathcal{L}^2(\Omega),\\
(c)&\;\int_{\Omega} |Du|^2\le \liminf_{k'\to\infty}\int_{\Omega}|Du_k'|^2.
\end{split}
\end{equation*}
\end{theorem}
\begin{proof}
There is actually nothing to say, once we notice that the convergence of $u_i$ is actually just equivalent
to the convergence of $\xi\circ u_i$.
\end{proof}
\begin{theorem}
If $\{u_j\}$ is a sequence of energy minimizing maps in
$\mathcal{Y}_2(\Omega, \underline{Q}(N))$ with
$\sup_{j}||u_j||_{\mathcal{Y}_2(B_\rho^m(y))}<\infty$ for each
ball $B_\rho^m(y)$ with $B_\rho^m(y) \subset\Omega$, then there is
a subsequence $\{u_{j'}\}$ and a minimizing map $u\in
\mathcal{Y}_2(\Omega,\underline{Q}(N))$ such that $u_{j'}\to u$ in
$\mathcal{Y}_2(B_\rho^m(y),\mathbb{Q})$ on each ball
$B_\rho^m(y)\subset \Omega$.
\end{theorem}
\begin{proof}
The proof is basically the same as in $\cite{sl},\S 2.9$.
\end{proof}
\section{Step minimizing sequences}
We assume that $f_0\in\mathcal{Y}_2(B_1^m,\mathbb{Q})$, $f$ is
strictly defined, and
$$Dir(f_0;B_1^m)<\infty, ||\xi\circ f_0||_{\mathcal{L}^2}<\infty$$
Define
$$\mathcal{M}=\{u\in\mathcal{Y}_2(B_1^m,\mathbb{Q}),\mbox{strictly defined}, u=f_0\;\mbox{on}\;\partial B_1^m\}$$
Given an $h>0$, we define inductively, for
$k=1,2,3,\cdot\cdot\cdot,$ a multiple-valued sobolev function
$f_h^k$ such that $f_h^0=f_0$ and $f_h^k$ minimizes the functional
(step minimizing with respect to $h$)
$$G(g)=G(f_h^{k-1},g,h)=
Dir(g;B_1^m)+\DF{1}{h/2^k}||\xi\circ f_h^{k-1}-\xi\circ
g||_{\mathcal{L}^2}^2$$ where $g\in\mathcal{M}$.
\begin{proposition}$f_h^k$ is defined for each
positive integer $k$ and each $h>0$.
\end{proposition}
\begin{proof}
Let
$$L_k=\inf\{G(f_h^{k-1},g,h):g\in\mathcal{M}\},k=1,2,3,\cdot\cdot\cdot$$
Since $G(f_h^0,f_0,h)=Dir(f_0;B_1^m)<\infty,$ $L_1<\infty.$\\
Let $\{g_i\}\subset \mathcal{M}$ be a sequence such that
$$\lim_{i\rightarrow \infty}G(f_h^0,g_i,h)=L_1,$$
then $\sup_i G(f_h^0,g_i,h)<\infty$.\\
It is easy to see
$$Dir(g_i;B_1^m)\le G(f_h^0,g_i,h)\le \sup_i G(f_h^0,g_i,h)<\infty,$$
By Theorem 2.5, there is a subsequence, still denoted as $g_i$
such that
$$g_i \rightharpoonup g\in \mathcal{Y}_2,$$
which gives
\begin{equation*}
\begin{split}
\int_{B_1^m} |\xi\circ f_h^0-\xi\circ g|^2 dx&=\int_{B_1^m} \lim_{i\rightarrow \infty} |
\xi\circ f_h^0-\xi\circ g_i|^2 dx\\
&\le \liminf_{i\rightarrow \infty} \int_{B_1^m}|\xi\circ
f_h^0-\xi\circ g_i|^2 dx\;(\mbox{by Fatou lemma})
\end{split}
\end{equation*}
and
$$Dir(g;B_1^m)\le \liminf_{i\to \infty} Dir(g_i;B_1^m)<\infty$$
with $g=f_h^0\;\mbox{on}\;\partial B_1^m$. By re-define $g$ at
some points if necessary, we may assume $g\in\mathcal{M}$.
Therefore
\begin{equation*}
\begin{split}
G(f_h^0,g,h)&=Dir(g;B_1^m)+\DF{\int_{B_1^m}|\xi\circ
f_h^0-\xi\circ
g|^2 dx}{h/2}\\
&\le \liminf_{i\to \infty} Dir(g_i;B_1^m)+\DF{\liminf_{i\to
\infty} \int_{B_1^m}|\xi\circ
f_h^0-\xi\circ g_i|^2 dx}{h/2}\\
&\le \liminf_{i\to\infty} G(f_h^0,g_i,h)=L_1
\end{split}
\end{equation*}
Hence $G(f_h^0,g,h)=L_1$. We just let $f_h^1$ to be $g$.\\
Now we assume that $f_h^k$ exists and $Dir(f_h^{k};B_1^m)<\infty.$
$$L_{k+1}=\inf\{G(f_h^k,g,h):g\in\mathcal{M}\}\le G(f_h^k,f_h^k,h)=Dir(f_h^k;B_1^m)<\infty.$$
Let $\{g_i\}\subset \mathcal{M}$ be a sequence that
$$\lim_{i\rightarrow \infty}G(f_h^k,g_i,h)=L_{k+1},$$
then $\sup_i G(f_h^k,g_i,h)<\infty$.\\
It is easy to see
$$Dir(g_i;B_1^m)\le G(f_h^k,g_i,h)\le \sup_i G(f_h^k,g_i,h)<\infty,$$
By Theorem 2.5, there is a subsequence, still denoted as $g_i$
such that
$$g_i \rightharpoonup g\in \mathcal{Y}_2,$$
which gives
\begin{equation*}
\begin{split}
\int_{B_1^m} |\xi\circ f_h^k-\xi\circ g|^2 dx&=\int_{B_1^m}
\lim_{i\rightarrow \infty} |
\xi\circ f_h^k-\xi\circ g_i|^2 dx\\
&\le \liminf_{i\rightarrow \infty} \int_{B_1^m}|\xi\circ
f_h^k-\xi\circ g_i|^2 dx\;(\mbox{by Fatou lemma})
\end{split}
\end{equation*}
and
$$Dir(g;B_1^m)\le \liminf_{i\to \infty} Dir(g_i;B_1^m)<\infty$$
with $g=f_h^0\;\mbox{on}\;\partial B_1^m$. By re-define $g$ at
some points if necessary, we may assume $g\in\mathcal{M}$.
Therefore
\begin{equation*}
\begin{split}
G(f_h^k,g,h)&=Dir(g;B_1^m)+\DF{\int_{B_1^m}|\xi\circ
f_h^k-\xi\circ
g|^2 dx}{h/2^{k+1}}\\
&\le \liminf_{i\to \infty} Dir(g_i;B_1^m)+\DF{\liminf_{i\to
\infty} \int_{B_1^m}|\xi\circ
f_h^k-\xi\circ g_i|^2 dx}{h/2^{k+1}}\\
&\le \liminf_{i\to\infty} G(f_h^k,g_i,h)=L_{k+1}
\end{split}
\end{equation*}
Hence $G(f_h^k,g,h)=L_{k+1}$. We just let $f_h^{k+1}$ to be $g$.\\
Repeat the process, we will prove the proposition.
\end{proof}
\begin{proposition}
The energy of the sequence $\{f_h^k\}$ is non-increasing,
$$Dir(f_h^k;B_1^m)\le Dir(f_h^{k-1};B_1^m)\le Dir(f_0;B_1^m).$$
We also have the following estimate:
$$||\xi\circ f_h^{k-1}-\xi\circ f_h^k||_{\mathcal{L}^2}^2\le \DF{h}{2^k}(Dir(f_h^{k-1};B_1^m)-Dir(f_h^k;B_1^m)).$$
\end{proposition}
\begin{proof}
The proof is quite straightforward once we notice that
$$G(f_h^{k-1},f_h^k,h)\le G(f_h^{k-1},f_h^{k-1},h)=Dir(f_h^{k-1};B_1^m).$$
\end{proof}
\section{Construction of the flow}
Fix $h>0$, we will construct a multiple-valued function $F_h$ on the
cylinder $[0,\infty)\times B_1^m$ such that
$F_h(0,x)=f_0(x), x\in B_1^m$, and $F_h(t,x)=f_0(x), t\in [0,\infty), x\in \partial B_1^m$.\\
When $t\in[(i-1)h,ih]$, $i=1,2,3,\cdot\cdot\cdot$,
$$F_h(t,x):=f_h^{i-1}(x), \;\mbox{if}\;t\in[(i-1)h,ih-\DF{h}{2^i}]$$
$$F_h(t,x):=\xi^{-1}\circ \rho\circ
[\DF{ih-t}{h/2^i}\xi\circ
f_h^{i-1}+\DF{t-(ih-\DF{h}{2^i})}{h/2^i}\xi\circ
f_h^i]\;\mbox{if}\;t\in[ih-\DF{h}{2^i},ih]$$ It is easy to see that
the function $F_h$ is well-defined by using the fact that
$\rho\circ\xi=\xi$.\\
As for the boundary date of $F_h$, we only have to check the
boundary date on $[0,\infty)\times \partial B_1^m$. Take $t\in
[ih-h/2^i,ih], x\in\partial B_1^m$, since
$f_h^{i-1}(x)=f_h^i(x)=f_0(x),$
\begin{equation*}
\begin{split}
F_h(t,x)&=\xi^{-1}\circ \rho\circ [\DF{ih-t}{h/2^i}\xi\circ
f_0(x)+\DF{t-(ih-\DF{h}{2^i})}{h/2^i}\xi\circ f_0(x)]\\
&=\xi^{-1}\circ\rho \circ \xi\circ f_0(x)=f_0(x).
\end{split}
\end{equation*}
Obviously, $F_h\in\mathcal{Y}_2([0,\infty)\times B_1^m,\mathbb{Q})$. \\
Denoting by $F_k$ the function $F_{1/2^k}$. We will show that for any $T>0$,
$$\sup_k Dir(F_k;[0,T]\times B_1^m)<\infty.$$
Choose positive integer $N$ such that $(N-1)h\le T<Nh$, where $h=1/2^k$.\\
Since $F_h(t,x)=f_h^{i-1}(x),t\in[(i-1)h,ih-\DF{h}{2^i}]$,
\begin{equation*}
\begin{split}
\mbox{Energy of}\;F_h\;\mbox{over}\;[(i-1)h,ih-\DF{h}{2^i}]\times B_1^m&=(h-\DF{h}{2^i})\;
\mbox{Energy of}\;f_h^{i-1}\\
&\le h\;\mbox{Energy of}\;f_h^{i-1}\\
&\le h\;\mbox{Energy of}\; f_0
\end{split}
\end{equation*}
Sum them up, we get
$$\mbox{Energy of}\;F_h\;\mbox{over}\;
(\cup_i [(i-1)h,ih-\DF{h}{2^i}]\times B_1^m) \cap [0,T]\times B_1^m\le T\cdot\;\mbox{Energy of}\;f_0$$
As for the case when $t\in[ih-\DF{h}{2^i},ih],i=1,2,3,\cdot\cdot\cdot$,
\begin{equation*}
\begin{split}
|\DF{\partial (\xi\circ F_h)}{\partial t}|^2&\le (\mbox{Lip}\rho)^2|\DF{-1}{h/2^i}\xi\circ f_h^{i-1}+\DF{1}
{h/2^i}\xi\circ f_h^i|^2\\
&=\DF{(\mbox{Lip}\rho)^2}{(h/2^i)^2}|\xi\circ f_h^{i-1}-\xi\circ f_h^i|^2.
\end{split}
\end{equation*}
Integrating of that gives
\begin{equation*}
\begin{split}
\int_{ih-h/2^i}^{ih} dt\int_{B_1^m} |\DF{\partial (\xi\circ F_h)}{\partial t}|^2 dx & \le
\int_{ih-h/2^i}^{ih} dt\int_{B_1^m} \DF{(\mbox{Lip}\rho)^2}{(h/2^i)^2}|\xi\circ f_h^{i-1}-\xi\circ f_h^i|^2 dx\\
&=(\mbox{Lip}\rho)^2\DF{\int_{B_1^m} |\xi\circ f_h^{i-1}-\xi\circ f_h^i|^2 dx}{h/2^i}\\
&\le (\mbox{Lip}\rho)^2 (Dir(f_h^{i-1};B_1^m)-Dir(f_h^i;B_1^m))
\end{split}
\end{equation*}
As for $\DF{\partial (\xi\circ F_h)}{\partial x}$, we have
\begin{equation*}
\begin{split}
|\DF{\partial (\xi\circ F_h)}{\partial x}|^2&\le \;(\mbox{Lip}\rho)^2 |\DF{ih-t}{h/2^i}\DF{\partial (\xi\circ f_h^{i-1})}{\partial x}+
\DF{t-(ih-h/2^i)}{h/2^i}\DF{\partial (\xi\circ f_h^i)}{\partial x}|^2\\
&\le 2(\mbox{Lip}\rho)^2[(\DF{ih-t}{h/2^i})^2|\DF{\partial (\xi\circ f_h^{i-1})}{\partial x}|^2+
(\DF{t-(ih-h/2^i)}{h/2^i})^2|\DF{\partial (\xi\circ f_h^i)}{\partial x}|^2]
\end{split}
\end{equation*}
Integrating of that gives
$$\int_{ih-h/2^i}^{ih} dt\int_{B_1^m} |\DF{\partial (\xi\circ F_h)}{\partial x}|^2 dx$$
\begin{equation*}
\begin{split}
&\le
2(\mbox{Lip}\rho)^2\int_{ih-\DF{h}{2^i}}^{ih} dt\int_{B_1^m} [
(\DF{ih-t}{h/2^i})^2|\DF{\partial (\xi\circ f_h^{i-1})}{\partial x}|^2+
(\DF{t-(ih-\DF{h}{2^i})}{h/2^i})^2|\DF{\partial (\xi\circ f_h^i)}{\partial x}|^2]
dx\\
&=2(\mbox{Lip}\rho)^2(\DF{h}{3\times 2^i})\int_{B_1^m} |\DF{\partial (\xi\circ f_h^{i-1})}{\partial x}|^2+
|\DF{\partial (\xi\circ f_h^i)}{\partial x}|^2 dx\\
&=2(\mbox{Lip}\rho)^2(\DF{h}{3\times 2^i})(Dir(f_h^{i-1};B_1^m)+Dir(f_h^i;B_1^m))\\
&\le 2(\mbox{Lip}\rho)^2(\DF{h}{3\times 2^i}) 2 Dir(f_0;B_1^m)
\end{split}
\end{equation*}
Hence
$$\mbox{Energy of}\;F_h\;\mbox{over}\;(\cup_i [ih-h/2^i,ih]\times B_1^m)\cap [0,T]\times B_1^m$$
\begin{equation*}
\begin{split}
&\le\sum_{i=1}^N [(\mbox{Lip}\rho)^2
(Dir(f_h^{i-1};B_1^m)-Dir(f_h^i;B_1^m))+2(\mbox{Lip}\rho)^2
(\DF{h}{3\times 2^i}) 2Dir(f_0;B_1^m)]\\
&\le (\mbox{Lip}\rho)^2 Dir(f_0;B_1^m)+\DF{4}{3}(\mbox{Lip}\rho)^2 Dir(f_0;B_1^m)(T+h)\\
&\le (\mbox{Lip}\rho)^2 Dir(f_0;B_1^m)+\DF{4}{3}(\mbox{Lip}\rho)^2
Dir(f_0;B_1^m)(T+1)< \infty
\end{split}
\end{equation*}
In summary, we have
$$Dir(F_h;[0,T]\times B_1^m)<\infty\;\mbox{uniformly for}\;h=1/2^k.$$
Using Theorem 2.5, we have \begin{theorem} There exists a
subsequence of $F_k$ converging weakly in $\mathcal{Y}_2$ to a
multiple-valued function $F\in\mathcal{Y}_2([0,\infty)\times
B_1^m,\mathbb{Q})$ such that
$$F(0,x)=f_0(x),x\in B_1^m,$$
$$F(t,x)=f_0(x),t\in [0,\infty), x\in \partial B_1^m.$$
\end{theorem}
\begin{definition}
Denote
$$F_t(\cdot):B_1^m\rightarrow \mathbb{Q}, \;F_t(x)=F(t,x),$$
for any $t\in[0,\infty)$.
\end{definition}
\begin{theorem}
For $\mathcal{L}^1$ almost every $t>0$,
$$Dir(F_t;B_1^m)\le Dir(f_0;B_1^m).$$
\end{theorem}
\begin{theorem}
The flow is H$\ddot{o}$lder continuous with respect to the
$\mathcal{L}^2$ norm, i.e.,
$$||\xi\circ F_t-\xi\circ F_s||_{\mathcal{L}^2}\le \sqrt{s-t}\sqrt{Dir(f_0;B_1^m)}$$
for $0\le t<s$ for $\mathcal{L}^1$ almost every $t,s$.
\end{theorem}
\section{Proof of main theorems}
\begin{lemma}Suppose $A_i,i=1,2,3,\cdot\cdot\cdot$ are measurable sets in
$\mathbb{R}$, $\sum_{i=1}^\infty \mathcal{L}^1 (A_i)<\infty$, then
$$\mathcal{L}^1 (\overline{\lim} A_k)=0.$$
\end{lemma}
\begin{proof}
This is a fundamental result whose proof can be found in almost
any real analysis book.
\end{proof}
\begin{lemma}Suppose $u_i,u\in\mathcal{Y}_2([0,\infty)\times
B_1^m,\mathbb{R}^{PQ})$, $u_i\rightharpoonup u
\;\mbox{in}\;\mathcal{Y}_2$, $\sup_i
||u_i||_{\mathcal{Y}_2}<\infty$. Then for $\mathcal{L}^1$ almost
every $t>0$, there is a subsequence $i_k$, such that
$u_{i_k}(t,\cdot),u(t,\cdot)\in\mathcal{Y}_2(B_1^m,\mathbb{R}^{PQ})$
and
$$u_{i_k}(t,\cdot)\rightharpoonup u(t,\cdot)\;\mbox{in}\;\mathcal{Y}_2(B_1^m,\mathbb{R}^{PQ}).$$
\end{lemma}
\begin{proof}
It suffices to prove in the case that $u=0$ and the domain is $[0,1]\times B_1^m$.\\
Suppose
$$||u_i||_{\mathcal{Y}_2}^2=\int_0^1 dt\int_{B_1^m} |u_i|^2+|\nabla u_i|^2 dx\le M<\infty$$
Since $|\nabla u_i(t,\cdot)|\le|\nabla u_i|$,
\begin{equation*}
\begin{split}
\int_0^1 \underline{\lim} ||u_i(t,\cdot)||_{\mathcal{Y}_2}^2 dt&\le \underline{\lim}\int_0^1 ||u_i(t,\cdot)||_{\mathcal{Y}_2}^2 dt\\
&\le \underline{\lim}\int_0^1 dt \int_{B_1^m} |u_i|^2+|\nabla u_i|^2
dx\le M<\infty.
\end{split}
\end{equation*}
Therefore, for $\mathcal{L}^1$ almost every $t>0$, $\underline{\lim}||u_i(t,\cdot)||_{\mathcal{Y}_2}^2<\infty$. By definition of $
\lim\inf$, for such $t$, there is a subsequence $i'$(which may depend on $t$), such that
$$\lim_{i'\rightarrow \infty}||u_{i'}(t,\cdot)||_{\mathcal{Y}_2}^2<\infty$$
By Theorem 4.2, there is a subsequence $i''$ such that
$$u_{i''}(t,\cdot)\rightharpoonup h_t,\;\mbox{for some}\;h_t\in\mathcal{Y}_2(B_1^m,\mathbb{R}^{PQ})$$
Hence $u_{i''}(t,\cdot)\rightarrow h_t\;\mbox{in}\;\mathcal{L}^2$. We will show that $h_t=0$.\\
\begin{equation*}
\begin{split}
\int_0^1 \underline{\lim}||u_{i''}(t,\cdot)||_{\mathcal{L}^2}^2 dt&\le \underline{\lim}
\int_0^1 ||u_{i''}(t,\cdot)||_{\mathcal{L}^2}^2 dt\\
&=\underline{\lim}||u_{i''}||_{\mathcal{L}^2}^2\rightarrow 0
\end{split}
\end{equation*}
where the last limit comes from the assumption that $u_{i''}\rightharpoonup 0\;\mbox{in}\;\mathcal{Y}_2$.\\
So for $\mathcal{L}^1$ almost every $t>0$, $\underline{\lim}||u_{i''}(t,\cdot)||_{\mathcal{L}^2}^2=0$.\\
For such $t$, there is a subsequence $i'''$ such that
$$u_{i'''}(t,\cdot)\rightarrow 0\;\mbox{in}\;\mathcal{L}^2.$$
This proves the lemma.
\end{proof}
\subsection{Proof of Theorem 6.2}
\begin{proof}
Define
$$A_k=\cup_{i=1}^\infty [ih-h/2^i,ih],h=1/2^k$$
$$\mathcal{L}^1(A_k)=\sum_{i=1}^\infty \DF{h}{2^i}=h=1/2^k$$
$$\sum_{k=1}^\infty \mathcal{L}^1(A_k)=\sum_{k=1}^\infty \DF{1}{2^k}=1$$
By Lemma 7.1, $\mathcal{L}^1(\overline{\lim}A_k)=0$. \\Apply Lemma
7.2 to $F_k$, there is a subset $B\subset[0,\infty)$, with
$\mathcal{L}^1(B)=0$, such that for any $t\notin B$, there is a
subsequence (still denoted as $k$) such that
$$F_k(t,\cdot)\rightharpoonup
F_t\;\mbox{in}\;\mathcal{Y}_2.$$
Noticing
$$(\overline{\lim}A_k)^c=\underline{\lim}(A_k^c)=\{t:\;\mbox{there exists}\;n_0,\;\mbox{when}\;k\ge n_0,x\in A_k^c\}$$
When $t\notin B\cup \overline{\lim}A_k$, after finite steps,
$t\notin A_k$ for any $k$, i.e., after finite steps,
$$F_k(t,x)=f_h^{l-1}(x), \;\mbox{for}\;h=1/2^k, (l-1)h\le t< lh$$
Therefore
$$Dir(F_t;B_1^m)\le
\underline{\lim}Dir(F_k(t,\cdot);B_1^m)=\underline{\lim}Dir(f_h^{l-1};B_1^m)\le
Dir(f_0;B_1^m)$$ for any $t\notin B\cup \overline{\lim}A_k$.
\end{proof}
\subsection{Proof of Theorem 6.3}
\begin{proof}
Take $t,s\notin B\cup\overline{\lim}A_k$, $t<s$, there is a subsequence (still denoted as $k$) such that
$$F_k(t,\cdot)\rightarrow F_t, F_k(s,\cdot)\rightarrow F_s \;\mbox{in}\;\mathcal{L}^2.$$
After some finite steps,
$$F_k(t,x)=f_h^{l-1}(x), \;\mbox{for}\;h=1/2^k, (l-1)h\le t< lh$$
$$F_k(s,x)=f_h^{l'-1}(x), \;\mbox{for}\;h=1/2^k, (l'-1)h\le s< l'h.$$
Therefore $||\xi\circ F_k(t,\cdot)-\xi\circ F_k(s,\cdot)||_{\mathcal{L}^2}^2=||\xi\circ
f_h^{l-1}-\xi\circ f_h^{l'-1}||_{\mathcal{L}^2}^2$ when $k$ is big enough.\\
Using the basic inequality
$$(\sum_{i=1}^N a_i)^2\le N\sum_{i=1}^N a_i^2$$
we have
\begin{equation*}
\begin{split}
||\xi\circ F_k(t,\cdot)-\xi\circ F_k(s,\cdot)||_{\mathcal{L}^2}^2&=||\xi\circ
f_h^{l-1}-\xi\circ f_h^{l'-1}||_{\mathcal{L}^2}^2\\
&\le(l'-l)\sum_{i=l-1}^{l'-2}||\xi\circ f_h^i-\xi\circ f_h^{i+1}||_{\mathcal{L}^2}^2\\
&\le(l'-l)\sum_{i=l-1}^{l'-2}\DF{h}{2^{i+1}}(Dir(f_h^i;B_1^m)-Dir(f_h^{i+1};B_1^m))\\
&\le(l'-l)h\sum_{i=l-1}^{l'-2}(Dir(f_h^i;B_1^m)-Dir(f_h^{i+1};B_1^m))\\
&\le (l'-l)h Dir(f_0;B_1^m)\le (s-t+h)Dir(f_0;B_1^m)
\end{split}
\end{equation*}
The theorem is proved once we let $k\rightarrow \infty$ in the above inequality.
\end{proof}
\section{Special cases of this flow}
\subsection{Further properties about minimizers $f_h^k$}
In this section, we will look at the case when $n=1$, namely,
$\mathbb{Q}=\mathbb{Q}_Q(\mathbb{R})$. In this case, we have, for
any $S=\sum_{i=1}^Q [[s_i]],W=\sum_{i=1}^Q [[w_i]]$ such that
$-\infty<s_1\le s_2\le\cdot\cdot\cdot\le s_Q<\infty,-\infty<w_1\le
w_2\le\cdot\cdot\cdot\le w_Q<\infty$,
$$|\xi(S)-\xi(W)|^2=\mathcal{G}^2(S,W)=\sum_{i=1}^Q |s_i-w_i|^2.$$
\begin{theorem}
The minimizer $f_h^k$ satisfies the following equation:
$$\int_{B_1^m}
<\phi,\DF{\eta(f_h^k)-\eta(f_h^{k-1})}{h/2^k}>+<D\phi,D(\eta(f_h^k))>
d\mathcal{L}^m=0,$$ for any $\phi\in C_0^1(B_1^m,\mathbb{R})$.
\end{theorem}
\begin{proof}
For simplicity, denote $$f_h^k(x)=\sum_{i=1}^Q
[[f_i(x)]],f_h^{k-1}(x)=\sum_{i=1}^Q [[f_i^0(x)]].$$ Take any
smooth function $\phi\in C_0^1(B_1^m,\mathbb{R})$, let
$$u_t(x)=f_h^k(x)+tQ[[\phi(x)]]=\sum_{i=1}^Q [[f_i(x)+t\phi(x)]]\in\mathcal{M}.$$
Using Theorem 2.6, we have
\begin{equation*}
\begin{split}
Dir(u_t;B_1^m)&=Dir(f_h^k;B_1^m)+2Q\int_{B_1^m} <D(\eta\circ
f_h^k),D(t\phi)>d\mathcal{L}^m+QDir(t\phi;B_1^m)\\
&=Dir(f_h^k;B_1^m)+2tQ\int_{B_1^m}<D(\eta\circ
f_h^k),D\phi>d\mathcal{L}^m+t^2QDir(\phi;B_1^m)
\end{split}
\end{equation*}
Fix any permutation $\sigma$ of $\{1,2,\cdot\cdot\cdot,Q\}$,
$$\sum_{i=1}^Q |f_i(x)+t\phi(x)-f_{\sigma(i)}^0(x)|^2$$
\begin{equation*}
\begin{split}
&=\sum_{i=1}^Q
|f_i(x)-f_{\sigma(i)}^0(x)|^2+t^2|\phi(x)|^2+2t<\phi(x),f_i(x)-f_{\sigma(i)}^0(x)>\\
&=\sum_{i=1}^Q
|f_i(x)-f_{\sigma(i)}^0(x)|^2+t^2Q|\phi(x)|^2+2t<\phi(x),\sum_{i=1}^Q
f_i(x)-\sum_{i=1}^Q f_{\sigma(i)}^0(x)>\\
&=\sum_{i=1}^Q
|f_i(x)-f_{\sigma(i)}^0(x)|^2+t^2Q|\phi(x)|^2+2t<\phi(x),Q(\eta\circ
f_h^k-\eta\circ f_h^{k-1})>
\end{split}
\end{equation*}
Therefore,
$$\mathcal{G}^2(u_t,f_h^{k-1})=\mathcal{G}^2(f_h^k,f_h^{k-1})+t^2Q|\phi|^2+2tQ<\phi,\eta\circ
f_h^k-\eta\circ f_h^{k-1}>$$ Hence
$$G(f_h^{k-1},u_t,h)=Dir(u_t;B_1^m)+\int_{B_1^m}\mathcal{G}^2(u_t,f_h^{k-1})d\mathcal{L}^m/(h/2^k)$$
$$=Dir(f_h^k;B_1^m)+2tQ\int_{B_1^m}<D(\eta\circ
f_h^k),D\phi>d\mathcal{L}^m+t^2QDir(\phi;B_1^m)+$$
$$[\int_{B_1^m}\mathcal{G}^2(f_h^k,f_h^{k-1})+t^2Q|\phi|^2+2tQ<\phi,\eta\circ
f_h^k-\eta\circ f_h^{k-1}>d\mathcal{L}^m]/(h/2^k)$$ Since
$u_0=f_h^k$ and $f_h^k$ minimizes the functional
$G(f_h^{k-1},g,h)$,
$$0=\DF{dG(f_h^{k-1},u_t,h)}{dt}|_{t=0}=$$
$$2Q\int_{B_1^m}<D(\eta\circ f_h^k),D\phi>d\mathcal{L}^m+2Q\int_{B_1^m} <\phi,\eta\circ f_h^k-
\eta\circ f_h^{k-1}>d\mathcal{L}^m/(h/2^k),$$ which means
$$\int_{B_1^m}
<\phi,\DF{\eta(f_h^k)-\eta(f_h^{k-1})}{h/2^k}>+<D\phi,D(\eta(f_h^k))>
d\mathcal{L}^m=0,$$ for any $\phi\in C_0^1(B_1^m,\mathbb{R})$.
\end{proof}
\subsection{The case when the initial data is symmetric}
In this section, $Q$ will be two, i.e,
$\mathbb{Q}=\mathbb{Q}_2(\mathbb{R})$.\\Given $g\in
C^{\infty}(B_1^m;\mathbb{R})$ such that $g$ is nonnegative. Define
$$f_0(x)=[[g(x)]]+[[-g(x)]]\in\mathcal{Y}_2(B_1^m;\mathbb{Q}_2(\mathbb{R}))$$
Now let's consider the flow with $f_0$ being the initial data.\\
Fix $h>0$, let
$$f_h^0=f_0.$$
Therefore $\eta(f_h^0)\equiv 0$. By Theorem 8.1, $\eta(f_h^1)$
satisfies the following equality:
$$\int_{B_1^m}
<\phi,\DF{\eta(f_h^1)}{h/2}>+<D\phi,D(\eta(f_h^1))>
d\mathcal{L}^m=0,$$ for any $\phi\in C_0^1(B_1^m,\mathbb{R})$.\\
That means $\eta(f_h^1)$ is a weak solution of this boundary-value
problem
$$\left\{
\begin{array}{cccc}
-\Delta u+\DF{u}{h/2}=0\;\;\;\mbox{in}\;B_1^m\\
u=0\;\;\;\;\;\;\;\;\;\;\;\;\;\;\;\;\;\;\;\;\;\mbox{on}\;\partial
B_1^m
\end{array}\right.$$
From $\cite{el}, \S 6.2$, there is a unique weak solution to it.
Since obviously zero is a solution to the above problem, we get
$$\eta(f_h^1)\equiv 0.$$
The same argument gives $\eta(f_h^k)\equiv
0,k=1,2,\cdot\cdot\cdot$. \\In spirit of $\eta(f_h^1)\equiv 0$, we
can write $f_h^1$ as
$$f_h^1=[[f(x)]]+[[-f(x)]],$$
where $f\ge 0$. \\
Since $\xi\circ f_h^1=(-f(x),f(x))$ and $\xi\circ
f_h^1\in\mathcal{Y}_2(B_1^m,\mathbb{R}^2)$,
$$f\in\mathcal{Y}_2(B_1^m,\mathbb{R}).$$
Take any nonnegative function $\phi\in C_0^1(B_1^m,\mathbb{R})$,
consider
$$f_t(x)=[[f(x)+t\phi(x)]]+[[-f(x)-t\phi(x)]]\in\mathcal{M}.$$
We have
\begin{equation*}
\begin{split}
Dir(f_t;B_1^m)&=2\int_{B_1^m}|Df+tD\phi|^2 dx\\
&=2\int_{B_1^m} [|Df|^2+t^2|D\phi|^2+2tDf\cdot D\phi]dx
\end{split}
\end{equation*}
\begin{equation*}
\begin{split}
\DF{1}{h/2}\int_{B_1^m} \mathcal{G}^2(f_h^0,f_t)
dx&=\DF{4}{h}\int_{B_1^m} |f+t\phi-g|^2dx\\
&=\DF{4}{h}\int_{B_1^m} [|f-g|^2+t^2\phi^2+2t\phi(f-g)]dx
\end{split}
\end{equation*}
Therefore,
$$0=\DF{d}{dt}|_{t=0}G(f_h^0,f_t,h)=4\int_{B_1^m} [Df\cdot D\phi+\DF{1}{h/2}\phi(f-g)]dx$$
for any nonnegative $\phi\in C_0^1(B_1^m,\mathbb{R})$. Because of
the linearity of the above equation, we conclude that $f$ is a
weak solution of the following boundary-value problem:
$$\left\{
\begin{array}{c}
-\Delta u+\DF{u}{h/2}=\DF{g}{h/2}\;\;\;\mbox{in}\;B_1^m\\
u=g\;\;\;\;\;\;\;\;\;\;\;\;\;\;\;\;\;\;\;\;\;\mbox{on}\;\partial
B_1^m
\end{array}\right.$$
By introducing $\tilde{u}=f-g$, we see $\tilde{u}$ is a weak
solution of the following boundary-value problem:
$$\left\{
\begin{array}{c}
-\Delta \tilde{u}+\DF{\tilde{u}}{h/2}=\Delta g\;\;\;\mbox{in}\;B_1^m\\
\tilde{u}=0\;\;\;\;\;\;\;\;\;\;\;\;\;\;\;\;\;\;\;\;\;\mbox{on}\;\partial
B_1^m
\end{array}\right.$$
which has a unique weak solution $\tilde{u}$. Moreover by the
regularity theorem in $\cite{el},\S 6.3$,
$$\tilde{u}\in C^{\infty}(B_1^m)$$
Hence so is $f$. Hence $f$ actually is a smooth solution of the
following PDE:
$$\DF{f-g}{h/2}=\Delta f.$$
Now let us denote:
$$f_h^k(x)=[[\tilde{f_h^k}(x)]]+[[-\tilde{f_h^k}(x)]],\tilde{f_h^k}(x)\ge 0, k=0,1,2,\cdot\cdot\cdot.$$
From the previous argument, we know each $\tilde{f_h^k}$ is a
smooth solution of the following PDE:
$$\DF{\tilde{f_h^k}-\tilde{f_h^{k-1}}}{h/2^k}=\Delta \tilde{f_h^k}$$
This gives the proof of this theorem:
\begin{theorem}
Suppose
$f_0=[[g]]+[[-g]]\in\mathcal{Y}_2(B_1^m,\mathbb{Q}_2(\mathbb{R}))$,
with nonnegative function $g\in C^{\infty}(B_1^m,\mathbb{R})$.
Then $\eta(f_h^k)\equiv 0$ for $k=0,1,2,\cdot\cdot\cdot.$
Moreover, if we denote
$$f_h^k(x)=[[\tilde{f_h^k}(x)]]+[[-\tilde{f_h^k}(x)]],\tilde{f_h^k}(x)\ge 0, k=0,1,2,\cdot\cdot\cdot.$$
then each $\tilde{f_h^k}\in C^{\infty}(B_1^m,\mathbb{R})$ and
satisfies the following PDEs:
$$\DF{\tilde{f_h^k}-\tilde{f_h^{k-1}}}{h/2^k}=\Delta \tilde{f_h^k}$$
\end{theorem}
Next, we will show that $f_h^k$ has no branch points. Namely,
\begin{theorem}
With the same assumptions as the above theorem, if moreover, $g$
is not identically zero, then
$$\tilde{f_h^k}(x)>0,x\in(B_1^m)^{\circ},k=1,2,\cdot\cdot\cdot.$$
In particular, $f_h^k(x)\not= 2[[0]]$, for any $x\in
(B_1^m)^{\circ}$, $k=1,2,\cdot\cdot\cdot$.
\end{theorem}
\begin{proof}
From Theorem 8.2, $\tilde{f_h^1}$ satisfies the following PDE:
$$\DF{\tilde{f_h^1}-\tilde{f_h^0}}{h/2}-\Delta \tilde{f_h^1}=0.$$
Consider this operator:
$$Lu=-\Delta u+\DF{u}{h/2}.$$
Let $\psi=\tilde{f_h^1}$. We have
$$L\psi=\DF{\tilde{f_h^0}}{h/2}\ge 0$$
By strong maximum principle of $\cite{el},\S 6.4$, we conclude
that if $\psi$ attains a nonpositive minimum over $B_1^m$ at an
interior point, then $\psi$ is constant within $B_1^m$. Since
$\psi\equiv 0$ on $\partial B_1^m$, either $\psi>0$ in
$(B_1^m)^{\circ}$ or $\psi\equiv 0 \;\mbox{in}\;B_1^m$. But if
$\psi\equiv 0$ in $B_1^m$, then $\tilde{f_h^0}=(h/2)L\psi\equiv
0$, which contradicts to our assumption about $g$. Therefore
$$\psi>0\;\mbox{in}\;(B_1^m)^{\circ},$$
which means $\tilde{f_h^1}(x)>0,x\in
(B_1^m)^{\circ}$.\\
Now we will use induction to prove the theorem. We assume that we
have already showed that:
$$\tilde{f_h^k}(x)>0,x\in(B_1^m)^{\circ}$$ We will show that
$$\tilde{f_h^{k+1}}(x)>0,x\in(B_1^m)^{\circ}.$$
From Theorem 8.2, $\tilde{f_h^{k+1}}$ satisfies the following PDE:
$$\DF{\tilde{f_h^{k+1}}-\tilde{f_h^k}}{h/2^{k+1}}-\Delta \tilde{f_h^{k+1}}=0.$$
Consider the operator:
$$Lu=-\Delta u+\DF{u}{h/2^{k+1}},$$
and let $$\psi=\tilde{f_h^{k+1}}.$$ Since
$$L\psi=\DF{\tilde{f_h^k}}{h/2^{k+1}}>0\;\mbox{in}\;(B_1^m)^{\circ},$$
the strong maximum principle tells either $\psi>0$ in
$(B_1^m)^{\circ}$ or $\psi\equiv 0$ in $B_1^m$.\\
$\psi$ can not be identically zero because otherwise $
\tilde{f_h^k}=(h/2^{k+1})L\psi\equiv 0$, a contradiction to our
induction assumption. Therefore, $\psi>0$ in $(B_1^m)^{\circ}$,
i.e,
$$\tilde{f_h^{k+1}}>0\;\mbox{in}\;(B_1^m)^{\circ}.$$
\end{proof}
\section{Heat flow of multiple-valued functions}
\subsection{Construction of heat flow for single-valued functions
by method of discrete Morse flow} We assume that $f_0\in
C^{\infty} (B_1^m,\mathbb{R})$, and
$$Dir(f_0;B_1^m)<\infty, ||f_0||_{\mathcal{L}^2}<\infty.$$
Define
$$\mathcal{W}=\{u\in\mathcal{Y}_2(B_1^m,\mathbb{R}),\mbox{strictly defined}, u=f_0\;\mbox{on}\;\partial B_1^m\}$$
Fix time $T$ and positive integer $N$, we let $h=T/N$. We define
inductively, for $k=0,1,2,\cdot\cdot\cdot,N$ a single-valued
sobolev function $f_h^k$ such that $f_h^0=f_0$ and $f_h^k$
minimizes the functional (step minimizing with respect to $h$)
$$G(g)=G(f_h^{k-1},g,h)=
Dir(g;B_1^m)+\DF{1}{h}||f_h^{k-1}-g||_{\mathcal{L}^2}^2$$ where
$g\in\mathcal{W}$.\\
Using the same method as in Proposition 5.1 and 5.2, we have:
\begin{proposition}$f_h^k$ is defined for each
positive integer $k$ and each $h=T/N>0$.
\end{proposition}
\begin{proposition}
The energy of the sequence $\{f_h^k\}$ is non-increasing,
$$Dir(f_h^k;B_1^m)\le Dir(f_h^{k-1};B_1^m)\le Dir(f_0;B_1^m).$$
We also have the following estimate:
$$||f_h^{k-1}-f_h^k||_{\mathcal{L}^2}^2\le h(Dir(f_h^{k-1};B_1^m)-Dir(f_h^k;B_1^m)).$$
\end{proposition}
We can also prove the following result in a similar way as in
Theorem 8.2, \begin{theorem} For any positive integer $N$, $f_h^i$
is smooth for $i=0,1,2,\cdot\cdot\cdot N$, and $h=T/N$. Moreover
they satisfy the following PDEs:
$$\DF{f_h^{i+1}-f_h^i}{h}=\Delta f_h^{i+1},i=0,1,2,\cdot\cdot\cdot N-1.$$
\end{theorem}
We define a function $F_N:[0,T]\times B_1^m\rightarrow
\mathbb{R}$ as follows:
$$F_N(t,x)=[\DF{(i+1)h-t}{h}f_h^i(x)+\DF{t-ih}{h}f_h^{i+1}(x)],\;\mbox{when}\;t\in[ih,(i+1)h],x\in B_1^m,$$
where $i=0,1,\cdot\cdot\cdot,N-1$.\\
When $t\in[ih,(i+1)h]$,
$$\DF{\partial F_N}{\partial t}=\DF{1}{h}(f_h^{i+1}-f_h^i),$$
$$D_x F_N=[\DF{(i+1)h-t}{h}D f_h^i+\DF{t-ih}{h}D f_h^{i+1}].$$
Therefore,
\begin{equation*}
\begin{split}
\int_{ih}^{(i+1)h}\int_{B_1^m}|\DF{\partial F_N}{\partial t}|^2
dtdx&=\DF{1}{h}\int_{B_1^m}(f_h^i-f_h^{i+1})^2 dx\\
&\le Dir(f_h^i)-Dir(f_h^{i+1}),
\end{split}
\end{equation*}
and
\begin{equation*}
\begin{split}
\int_{ih}^{(i+1)h}\int_{B_1^m}|D_x F_N|^2 dtdx&\le
2[\int_{ih}^{(i+1)h}(\DF{(i+1)h-t}{h})^2 dt\int_{B_1^m} |Df_h^i|^2
dx\\
&+\int_{ih}^{(i+1)h}(\DF{t-ih}{h})^2 dt\int_{B_1^m} |Df_h^{i+1}|^2
dx]\\
&=\DF{2h}{3}(Dir(f_h^i)+Dir(f_h^{i+1}))\le \DF{4h}{3}Dir(f_0)
\end{split}
\end{equation*}
Integration over the domain $[0,T]\times B_1^m$ gives
\begin{equation*}
\begin{split}
\mbox{Energy of}\;F_N &\le
\sum_{i=0}^{N-1}(Dir(f_h^i)-Dir(f_h^{i+1}))+\sum_{i=0}^{N-1}
\DF{4h}{3}Dir(f_0)\\
&\le Dir(f_0)+\DF{4T}{3} Dir(f_0)<\infty
\end{split}
\end{equation*}
In spirit of Theorem 2.5, we have the following theorem,
\begin{theorem}There is a subsequence of $F_N$, still denoted as $F_N$ for
simplicity, such that
$$F_N\rightharpoonup F\in \mathcal{Y}_2([0,T]\times
B_1^m,\mathbb{R})$$ for some single-valued sobolev function $F$.
Moreover,
$$F(0,x)=f_0(x),x\in B_1^m,$$
$$F(t,x)=f_0(x),t\in[0,T],x\in\partial B_1^m.$$
\end{theorem}
Next we will show that $F\in C^{\infty}([0,T]\times
B_1^m,\mathbb{R})$ and it solves the heat equation
$$\DF{\partial F}{\partial t}=\Delta_x F.$$
\begin{lemma}
For any fixed $h>0$, $\sup_{x\in B_1^m} |f_h^i(x)|$ is
non-increasing in $i$.
\end{lemma}
\begin{proof}
We put $\rho_{n-1}=\sup_{x\in B_1^m} |f_h^{n-1}(x)|$, and
$v_n=\max\{|f_h^n|^2-\rho_{n-1}^2,0\}$. We will show that
$v_n=0\;a.e.$ in $B_1^m$.\\
Let $\eta\in C_0^{\infty}(B_1^m)$ be a nonnegative function. It is
easy to see that for any sufficiently small $\epsilon>0$,
$$\DF{1}{h\epsilon}\int_{B_1^m}
|f_h^n-f_h^{n-1}|^2-|(1-\epsilon\eta)f_h^n-f_h^{n-1}|^2 dx$$
$$=\DF{1}{h}\int_{B_1^m}
\eta\{(1-\epsilon\eta)|f_h^n|^2-|f_h^{n-1}|^2+|f_h^n-f_h^{n-1}|^2\}dx:=I$$
and hence
\begin{equation*}
\begin{split}
0&\ge
\DF{G(f_h^{n-1},f_h^n,h)-G(f_h^{n-1},(1-\epsilon\eta)f_h^n,h)}{\epsilon}=\DF{Dir(f_h^n)-Dir((1-\epsilon\eta)f_h^n)}{\epsilon}+I\\
&=\DF{1}{\epsilon}\int_{B_1^m}\{
(2\epsilon\eta-\epsilon^2\eta^2)|Df_h^n|^2-\epsilon^2|D\eta|^2|f_h^n|^2+\epsilon(1-\epsilon\eta)D[(f_h^n)^2]\cdot
D\eta \}dx+I
\end{split}
\end{equation*}
Let $\epsilon\to 0$,
\begin{equation*}
\begin{split}
0&\ge \int_{B_1^m} \{2\eta|Df_h^n|^2+D[(f_h^n)^2]\cdot D\eta\}
dx+\lim_{\epsilon\to 0} I\\
&\ge \int_{B_1^m} D[(f_h^n)^2]\cdot D\eta dx+\DF{1}{h}\int_{B_1^m}
\eta(|f_h^n|^2-|f_h^{n-1}|^2) dx
\end{split}
\end{equation*}
Setting $\eta=v_n$, the last inequality gives us that
\begin{equation*}
\begin{split}
0&\ge \DF{1}{h}\int_{\{x:|f_h^n(x)|>\rho_{n-1}\}}
(|f_h^n|^2-\rho_{n-1}^2)(|f_h^n|^2-|f_h^{n-1}|^2)dx\\
&\ge \DF{1}{h}\int_{\{x:|f_h^n(x)|>\rho_{n-1}\}}
(|f_h^n|^2-\rho_{n-1}^2)^2dx,
\end{split}
\end{equation*}
which proves $|f_h^n|\le \rho_{n-1}$ a.e. in $B_1^m$.
\end{proof}
\begin{lemma}
$\{F_n\}$ is uniformly bounded in $[0,T]\times B_1^m$.
\end{lemma}
\begin{proof}
For any positive integer $n$, let $h=T/n$,
$$F_n(t,x)=\DF{(i+1)h-t}{h}f_h^i(x)+\DF{t-ih}{h}f_h^{i+1}(x),\;\mbox{when}\;ih\le
t\le (i+1)h.$$ Hence,
\begin{equation*}
\begin{split}
|F_n(t,x)|&\le
\DF{(i+1)h-t}{h}|f_h^i(x)|+\DF{t-ih}{h}|f_h^{i+1}(x)|\\
&\le \DF{(i+1)h-t}{h}\sup_{x\in B_1^m}
|f_h^i(x)|+\DF{t-ih}{h}\sup_{x\in B_1^m}|f_h^{i+1}(x)|\\
&\le \DF{(i+1)h-t}{h}\sup_{x\in B_1^m}
|f_h^0(x)|+\DF{t-ih}{h}\sup_{x\in B_1^m}|f_h^0(x)|\\
&= \sup_{x\in B_1^m} |f_0|<\infty,
\end{split}
\end{equation*}
where the third inequality comes from Lemma 9.1 and the last
inequality follows from the smoothness of $f_0$.
\end{proof}
Applying Theorem K in $\cite{kk}$ gives us
\begin{lemma}
There exist constants $C>0$ and $\mu>0$ such that
$$|F_n(t,x)-F_n(t',y)|\le C(|t-t'|^\mu+|x-y|^{2\mu}),\forall (t,x),
(t',y)\in [0,T]\times B_1^m,$$ i.e. $\{F_n\}$ is uniformly
equicontinuous.
\end{lemma}
Since each $f_h^i$ is smooth, each $F_n$ is smooth. By using the
Arzela-Ascoli theorem, we have
\begin{theorem}
There is a subsequence $\{F_{k_j}\}$ of $\{F_n\}$ which uniformly
converges to $F$, for $F:[0,T]\times B_1^m\to \mathbb{R}$ as in
Theorem 9.2. Therefore $F$ is continuous.
\end{theorem}
\begin{proof}
The Arzela-Ascoli theorem guarantees that there are a continuous
function $H:[0,T]\times B_1^m\to \mathbb{R}$ and a subsequence
$\{F_{k_j}\}$ such that $F_{k_j}$ uniformly converges to $H$.
Since the $\mathcal{Y}_2$ converges implies $\mathcal{L}^2$
convergence in some subsequence, we know $F=H$.
\end{proof}
\begin{remark}
We will also just use the original sequence $\{F_n\}$ without
referring to any subsequence for simplicity.
\end{remark}
\begin{theorem}
$F$ is smooth and satisfies the heat equation $$\DF{\partial
F}{\partial t}=\Delta_x F.$$ \end{theorem} \begin{proof} Summing
up the PDEs in Theorem 9.1 gives us that
$$f_h^j-f_h^0=h\sum_{i=1}^{j-1} \Delta f_h^{i+1},j=1,2,\cdot\cdot\cdot,n$$
Fix $n$, let $h=T/n$. For any fixed $t\in [0,T]$, $x\in B_1^m$,
choose integer $l$ such that $lh\le t<(l+1)h$, then
$$F_n(t,x)=\DF{(l+1)h-t}{h}f_h^l(x)+\DF{t-lh}{h}f_h^{l+1}(x),$$
$$\Delta_x F_n(t,x)=\DF{(l+1)h-t}{h}\Delta f_h^l(x)+\DF{t-lh}{h}\Delta
f_h^{l+1}(x).$$ Therefore,
\begin{equation*}
\begin{split}
F_n(t,x)-F_n(0,x)&=\DF{(l+1)h-t}{h}f_h^l(x)+\DF{t-lh}{h}f_h^{l+1}(x)-f_h^0(x)\\
&=\DF{(l+1)h-t}{h}(f_h^l(x)-f_h^0(x))+\DF{t-lh}{h}(f_h^{l+1}(x)-f_h^0(x))\\
&=\DF{(l+1)h-t}{h}(\sum_{i=0}^{l-1} \Delta
f_h^{i+1})h+\DF{t-lh}{h}(\sum_{i=0}^l \Delta f_h^{i+1})h\\
&=h\sum_{i=0}^{l-1}\Delta f_h^{i+1}+(t-lh)\Delta f_h^{l+1}
\end{split}
\end{equation*}
Choose $\phi\in C_0^\infty(B_1^m;\mathbb{R})$, do the integration
by parts twice, we have
$$\int_{B_1^m}(F_n(t,x)-F_n(0,x))\phi(x)dx=\sum_{i=0}^{l-1}
(\int_{B_1^m} f_h^{i+1}\Delta \phi dx)h+(t-lh)\int_{B_1^m}
f_h^{l+1}\Delta \phi dx$$ Letting $n\to \infty$ gives
$$\int_{B_1^m}
(F(t,x)-F(0,x))\phi(x)dx=\int_0^t\int_{B_1^m}F(s,x)\Delta
\phi(x)dxds,$$
$$\int_{B_1^m}\DF{\partial F}{\partial
t}\phi(x)dx=\int_{B_1^m}F(t,x)\Delta
\phi(x)dx=-\int_{B_1^m}\nabla_x F(t,x)\cdot \nabla \phi dx$$ In a
world, $F$ is a weak solution of the heat equation
$$\DF{\partial F}{\partial t}=\Delta_x F.$$
Since $F(0,x)=f_0$ is smooth, $F$ must be the unique smooth
solution of the heat equation:
$$\DF{\partial F}{\partial t}=\Delta_x F,(t,x)\in [0,T]\times B_1^m,$$
$$F(0,x)=f_0(x),x\in B_1^m,$$
$$F(t,x)=f_0(x),t\in[0,T],x\in\partial B_1^m.$$
\end{proof}
\subsection{Heat flow of multiple-valued functions} In this
section, we use a modified sequence of functional to construct a
flow for multiple-valued functions which is a generalization of
ordinary heat flow. \\
Suppose
$f_0=[[g]]+[[-g]]\in\mathcal{Y}_2(B_1^m,\mathbb{Q}_2(\mathbb{R}))$,
with nonnegative, nonconstant function $g\in
C^{\infty}(B_1^m,\mathbb{R})$. The functionals are
$$G(g)=G(f_h^{k-1},g,h)=
Dir(g;B_1^m)+\DF{1}{h}||\xi\circ f_h^{k-1}-\xi\circ
g||_{\mathcal{L}^2}^2$$ where $g\in\mathcal{M}$.\\
Using the same argument in Theorem 8.3, we know that the flow
instantly separates. If we abuse our notation by denoting
$f_h^k=[[f_h^k]]+[[-f_h^k]]$, for nonnegative function $f_h^k\in
\mathcal{Y}_2(B_1^m,\mathbb{R})$, then each $f_h^k$ minimizes the
functional
$$G(g)=G(f_h^{k-1},g,h)=
Dir(g;B_1^m)+\DF{1}{h}||f_h^{k-1}-
g||_{\mathcal{L}^2}^2$$ where $g\in\mathcal{W}$.\\
We define the discrete flow as follows:
\begin{equation*}
\begin{split}
F_n(t,x)&=[[\DF{(i+1)h-t}{h}f_h^i(x)+\DF{t-ih}{h}f_h^{i+1}(x)]]+\\
&[[-\DF{(i+1)h-t}{h}f_h^i(x)-\DF{t-ih}{h}f_h^{i+1}(x)]]
\end{split}
\end{equation*}
for $ih\le t<(i+1)h$, $h=T/n$.\\
The above section tells that the limit multiple-valued function
$F:[0,T]\times B_1^m\to \mathbb{Q}_2(\mathbb{R})$ is smooth, and
can be written as
$$F(t,x)=[[F_1(t,x)]]+[[-F_1(t,x)]]$$
for some nonnegative smooth function $F_1:[0,T]\times B_1^m \to
\mathbb{R}$ which satisfies the ordinary heat equation.

\end{document}